\documentclass[11pt]{article}

\usepackage{amssymb,latexsym}
\usepackage{amsmath}
\usepackage{cite}
\usepackage{eucal}
\usepackage[latin1]{inputenc}
\usepackage{bbm}
\usepackage{theorem}

\title{Morita Equivalence, Picard Groupoids and Noncommutative Field Theories}

\author{\textbf{Stefan Waldmann\thanks{Stefan.Waldmann@physik.uni-freiburg.de}} 
  \\[0.5cm]
  Fakultät für Mathematik und Physik\\
  Physikalisches Institut\\
  Albert-Ludwigs-Universität Freiburg\\
  Hermann Herder Strasse 3\\
  D-79104 Freiburg\\
  Germany
  }

\date{April 2003\\[0.5cm] FR-THEP 2003/06}

\renewcommand{\mathbb}[1]{\mathbbm{#1}}

\newcommand{\im}{{\mathrm i}}

\newcommand{\cc}[1]      {\overline{{#1}}}
\newcommand{\id}         {{\mathsf{id}}}
\newcommand{\tr}         {\mathop{{\mathsf{tr}}}}

\newcommand{\End}        {{\mathsf{End}}}
\newcommand{\ring}[1]    {{\mathsf{{#1}}}}

\newcommand{\Unit}       {\mathbb{1}}
\newcommand{\cl}         {\mathrm{cl}}
\newcommand{\tensor}     {\!\otimes\!}
\newcommand{\Pic}        {\mathrm{Pic}}
\newcommand{\Picc}       {\underline{\mathrm{Pic}}}

\newtheorem{lemma}{Lemma}
\newtheorem{proposition}[lemma]{Proposition}
\newtheorem{theorem}[lemma]{Theorem}

\newtheorem{definition}[lemma]{Definition}

\theorembodyfont{\rmfamily}

\newtheorem{remark}[lemma]{Remark}

\begin{document}

\maketitle

\begin{abstract}
In this article we review recent developments on Morita equivalence of
star products and their Picard groups. We point out the relations between
noncommutative field theories and deformed vector bundles which give
the Morita equivalence bimodules.
\end{abstract}

%
%

\section{Introduction: Noncommutative Field Theories}
\label{sec:intro}

Noncommutative field theory has recently become a very active field in
mathematical physics, see e.g.
\cite{doplicher.fredenhagen.roberts:1995a,jurco.schupp.wess:2000a,seiberg.witten:1999a,connes.douglas.schwarz:1998a,schomerus:1999a}
to mention just a few references. Many additional references can be
found in the recent review \cite{szabo:2001a:pre} as well as in these
proceedings.

The purpose of this note is to point out some mathematical structures
underlying the noncommutative field theories and the relations to
deformation quantization \cite{bayen.et.al:1978a} and Morita
equivalence of star products
\cite{bursztyn.waldmann:2001a,bursztyn.waldmann:2002a,jurco.schupp.wess:2002a}.

In the commutative framework the (matter) fields are geometrically
described by sections $\mathcal{E} = \Gamma^\infty(E)$ of some vector
bundle $E \to M$ over the space-time manifold $M$. We consider here a
complex vector bundle $E$. Since a field $\phi \in \mathcal{E}$ can be
multiplied by a function $f \in C^\infty(M)$ and since clearly $(\phi
f)g = \phi(fg)$ we obtain a (right) module structure of $\mathcal{E}$
over the algebra of smooth complex-valued functions $C^\infty(M)$.
Gauge transformations are encoded in the action of the sections of the
endomorphism bundle $\Gamma^\infty(\End(E))$, i.e. $A \in
\Gamma^\infty(\End(E))$ can be applied to a field $\phi$ by pointwise
multiplication $A\phi$. In this way, $\mathcal{E}$ becomes a
$\Gamma^\infty(\End(E))$ left module. Moreover, the action of
$\Gamma^\infty(\End(E))$ commutes with the action of $C^\infty(M)$
\begin{equation}
    \label{eq:bimodule}
    (A \phi) f = A (\phi f),
\end{equation}
whence the space of fields $\mathcal{E}$ becomes a \emph{bimodule}
over the algebras $\Gamma^\infty(\End(E))$ and $C^\infty(M)$.

In order to formulate not only the kinematics but also the dynamics we
need a Lagrange density $\mathcal{L}$ for $E$. Geometrically, this is
a function on the first jet bundle of $E$. A particular important
piece in $\mathcal{L}$ is the \emph{mass term} which is encoded in a
\emph{Hermitian fibre metric} $h_0$ for $E$. Recall that a Hermitian
fibre metric is a map
\begin{equation}
    \label{eq:HermitianMetric}
    h_0 : \mathcal{E} \times \mathcal{E} \to C^\infty(M)
\end{equation}
such that $h_0$ is $C^\infty(M)$-linear in the second argument,
$h_0(\phi,\psi) = \cc{h_0(\psi,\phi)}$ and one has the positivity 
\begin{equation}
    \label{eq:Hnullpos}
    h_0(\phi,\phi)(x) > 0 
    \quad
    \text{iff} 
    \quad
    \phi(x) \ne 0.
\end{equation}
Then the mass term in $\mathcal{L}$ is just $h_0(\phi,\phi)$ and the
last condition \eqref{eq:Hnullpos} is the positivity of the masses.
Note that such a Hermitian fibre metric is also used to encode
geometrically some polynomial interaction terms like $\phi^4$. Hence
it is of major importance to have a definiteness like
\eqref{eq:Hnullpos}.

We presented this well-known geometrical formulation, see e.g. the
textbook \cite{thirring:1990a}, in order to motivate now the
noncommutative analogs. The main idea is that at some scale (Planck,
etc.) the space-time itself behaves in a noncommutative fashion. One
way to encode this noncommutative nature is to consider a star product
$\star$ on $M$ which makes the algebra of functions $C^\infty(M)$ into
a noncommutative algebra. Here we consider formal star products for
convenience, see \cite{bayen.et.al:1978a} as well as
\cite{dito.sternheimer:2002a,gutt:2000a} for recent reviews and
further references.

Thus let $\pi$ be a Poisson tensor on the space-time $M$ and let
$\star$ be a formal star product for $\pi$, i.e. a
$\mathbb{C}[[\lambda]]$-bilinear associative multiplication for
$C^\infty(M)[[\lambda]]$,
\begin{equation}
    \label{eq:starprod}
    f \star g = \sum_{r=0}^\infty \lambda^r C_r(f,g),
\end{equation}
with some bidifferential operators $C_r$ such that $C_0(f,g) = fg$ is
the undeformed product and $C_1(f,g) - C_1(g,f) = \im\{f,g\}$ gives
the Poisson bracket corresponding to $\pi$. Moreover, we assume $f
\star 1 = f = 1 \star f$ and $\cc{f \star g} = \cc{g} \star \cc{f}$.
The formal parameter $\lambda$ corresponds to the scale where the
noncommutativity becomes important. Two star products are called
\emph{equivalent} if there is a formal series $T = \id +
\sum_{r=1}^\infty \lambda^r T_r$ of differential operators $T_r$ such
that $T(f \star g) = Tf \star' Tg$. See
\cite{dewilde.lecomte:1983b,fedosov:1996a,omori.maeda.yoshioka:1991a,kontsevich:1997:pre}
for existence and
\cite{nest.tsygan:1995a,bertelson.cahen.gutt:1997a,weinstein.xu:1998a,deligne:1995a,kontsevich:1997:pre}
for the classification of such star products up to equivalence.

In order to give a geometrical framework of noncommutative field
theories we want a deformed picture of the above bimodule structure.
Thus we look for a right module structure $\bullet$ on the space
$\Gamma^\infty(E)[[\lambda]]$ with respect to the algebra
$C^\infty(M)[[\lambda]]$. Thus $\bullet$ is a
$\mathbb{C}[[\lambda]]$-bilinear map
\begin{equation}
    \label{eq:BulletDef}
    \phi \bullet f = \sum_{r=0}^\infty \lambda^r R_r(\phi, f),
\end{equation}
where $R_r: \Gamma^\infty(E) \times C^\infty(M) \to \Gamma^\infty(E)$
is a bidifferential operator with $R_0(\phi,f) = \phi f$ and
\begin{equation}
    \label{eq:defbimod}
    (\phi \bullet f) \bullet g = \phi \bullet (f \star g)
    \quad
    \textrm{and}
    \quad
    \phi \bullet 1 = \phi.
\end{equation}
This gives the right module structure. But we also need an associative
deformation $\star'$ of $\Gamma^\infty(\End(E))$ and a left module
structure $\bullet'$ such that we have
\begin{equation}
    \label{eq:starprime}
    (A \star' B) \bullet' \phi = A \bullet' (B \bullet' \phi),
    \quad
    \Unit \bullet' \phi = \phi
    \quad
    \textrm{and}
    \quad
    A \bullet' (\phi \bullet f) = (A \bullet' \phi) \bullet f.
\end{equation}
This gives then a deformed bimodule structure on $\Gamma^\infty(E)$.

If we are interested in the analog of the Hermitian metric $h_0$
then we want a $\mathbb{C}[[\lambda]]$-sesquilinear deformation
$\boldsymbol{h} = \sum_{r=0}^\infty \lambda^r h_r$ of $h_0$ where
$h_r: \Gamma^\infty(E) \times \Gamma^\infty(E) \to C^\infty(M)$ such
that
\begin{gather}
    \label{eq:hbulletstar}
    \boldsymbol{h}(\phi, \psi \bullet f) 
    = \boldsymbol{h}(\phi, \psi) \star f, 
    \\
    \label{eq:hcc}
    \boldsymbol{h} (\phi, \psi) = \cc{\boldsymbol{h}(\psi,\phi)}, 
    \\
    \label{eq:hpos}
    \boldsymbol{h}(\phi,\phi) \; \textrm{is positive},
    \\
    \label{eq:AbullethHermitian}
    \boldsymbol{h}(A \bullet' \phi, \psi) 
    = \boldsymbol{h}(\phi, A^* \bullet' \psi).
\end{gather}
The positivitiy in \eqref{eq:hpos} is understood in the sense of
$^*$-algebras over ordered rings, see
\cite{bursztyn.waldmann:2001a}. In the case of a vector bundle this
just means that $\boldsymbol{h}(\phi,\phi)$ can be written as a sum of
squares $\sum_i \cc{f}_i \star f_i$.

Having this structure one obtains a \emph{framework} for
noncommutative field theories \emph{beyond} the usual formulations on
a flat space-time with trivial vector bundle, very much in the spirit
of Connes' noncommutative geometry \cite{connes:1994a}. To formulate a
physical theory one needs of course much more, like an action
principle, convergence in the deformation parameter $\lambda$, a
quantization of this still classical theory, etc. All these questions
shall not be addressed in this work. Instead, we shall focus on the
question whether and how one can prove existence, construct, and
classify the structures $\bullet$, $\bullet'$, $\star'$ and
$\boldsymbol{h}$ out of the given classical data and a given star
product $\star$. The case of a line bundle $E = L$ plays a
particularly interesting role as this corresponds exactly to (complex)
scalar fields.

%
%

\section{Deformation of Projective Modules}
\label{sec:deform}

There are several different ways to construct deformed versions of a
Hermitian vector bundle. We shall focus on a rather general algebraic
construction before discussing the other possibilities. Fundamental is
the well-known \emph{Serre-Swan theorem} \cite{swan:1962a} in its
smooth version: The $C^\infty(M)$-module of sections
$\Gamma^\infty(E)$ is a finitely generated projective module.
Moreover, the $C^\infty(M)$-linear module endomorphisms are just the
sections of the endomorphism bundle, i.e. $\Gamma^\infty(\End(E)) =
\End_{C^\infty(M)}(\Gamma^\infty(E))$. Hence one finds a projection
$P_0 = P_0^2 \in M_N(C^\infty(M))$ where $N$ is sufficiently large,
such that
\begin{equation}
    \label{eq:Eprojective}
    \Gamma^\infty(E) \cong P_0 C^\infty(M)^N
    \quad
    \textrm{and}
    \quad
    \Gamma^\infty(\End(E)) \cong P_0 M_N(C^\infty(M)) P_0.
\end{equation}
If $E$ is equipped with a Hermitian fiber metric $h_0$ then one can
even find a Hermititan projection $P_0 = P_0^2 = P_0^*$ such that with
the identification of \eqref{eq:Eprojective} the Hermitian fiber
metric becomes
\begin{equation}
    \label{eq:hnullnice}
    h_0 (\phi,\psi) = \sum_{i=1}^N \cc{\phi}_i\psi_i,
\end{equation}
where $\phi = (\phi_1, \ldots, \phi_N)$, $\psi = (\psi_1, \ldots,
\psi_N) \in C^\infty(M)^N$ are elements in $P_0 C^\infty(M)^N$,
i.e. they satisfy $P_0 \phi = \phi$, $P_0\psi = \psi$.

It is worth to look at this situation in general. Thus let
$\mathcal{A}$ be an associative algebra over a ring $\ring{C}$ and let
$\star$ be an associative formal deformation of $\mathcal{A}$. We
denote the deformed algebra by $\boldsymbol{\mathcal{A}} =
(\mathcal{A}[[\lambda]], \star)$. Now let $\mathcal{E}$ be a finitely
generated projective right module over $\mathcal{A}$ and let
$\End_{\mathcal{A}}(\mathcal{E})$ denote the $\mathcal{A}$-linear
endomorphisms of $\mathcal{E}$. Then one has the following result, see
\cite{bursztyn.waldmann:2000b}:
\begin{theorem}
    \label{theorem:DefProj}
    There exists a deformation $\bullet$ of $\mathcal{E}$ into a
    $\boldsymbol{\mathcal{A}}$-right module $(\mathcal{E}[[\lambda]],
    \bullet)$ which is unique up to equivalence such that
    $(\mathcal{E}[[\lambda]], \bullet)$ is finitely generated and
    projective over $\boldsymbol{\mathcal{A}}$ and
    $\End_{\boldsymbol{\mathcal{A}}}(\mathcal{E}[[\lambda]], \bullet)$
    is isomorphic as $\ring{C}[[\lambda]]$-module to
    $\End_{\mathcal{A}}(\mathcal{E})[[\lambda]]$.
\end{theorem}
Equivalence of two deformations $\bullet$ and $\tilde{\bullet}$ means
that there is a map $T = \id + \sum_{r=1}^\infty \lambda^r T_r$ with
$T(\phi \bullet a) = T(\phi) \tilde{\bullet} a$ for all $\phi \in
\mathcal{E}[[\lambda]]$ and $a \in \mathcal{A}[[\lambda]]$.

The idea of the proof consists in first deforming the projection $P_0$
into a projection $\boldsymbol{P}$ with respect to the deformed
product $\star$ by using the formula \cite[Eq. (6.1.4)]{fedosov:1996a}
\begin{equation}
    \label{eq:defprojformula}
    \boldsymbol{P} = \frac{1}{2} + \left(P_0 - \frac{1}{2}\right)
    \star \frac{1}{\sqrt[\star]{1 + 4 (P_0 \star P_0 - P_0)}}.
\end{equation}
Then the $\boldsymbol{\mathcal{A}}$-right module $\boldsymbol{P} \star
\boldsymbol{\mathcal{A}}^N$ is obviously a finitely generated and
projective $\boldsymbol{\mathcal{A}}$-module and it turns out that it
is isomorphic to $\mathcal{E}[[\lambda]]$ as
$\ring{C}[[\lambda]]$-module. Then the uniqueness of the deformation
$\bullet$ up to equivalence follows from the fact the $\boldsymbol{P}
\star \boldsymbol{\mathcal{A}}^N$ is projective again. Indeed, let
$\mathcal{E}$ be endowed with the trivial $\boldsymbol{\mathcal{A}}$
right module structure given by $\phi \cdot a = \phi a_0$ for $\phi
\in \mathcal{E}$ and $a = \sum_{r=0}^\infty \lambda^r a_r \in
\mathcal{A}[[\lambda]]$. Then the classical limit map $\cl:
(\mathcal{E}[[\lambda]], \bullet) \longrightarrow (\mathcal{E},\cdot)$
(setting $\lambda = 0$) is a module morphism for
$\boldsymbol{\mathcal{A}}$ right modules. The same holds for any other
deformation $\tilde{\bullet}$. Since the deformation $\bullet$ is
projective and since $\cl$ is obviously surjective this means we can
find a module morphism $T:(\mathcal{E}[[\lambda]], \bullet)
\longrightarrow (\mathcal{E}[[\lambda]], \tilde{\bullet})$ such that
$\cl \circ T = \cl$. This implies $T = \id + \sum_{r=1}^\infty
\lambda^r T_r$ whence we have found an equivalence, see
\cite{bursztyn.waldmann:2000b} for details.

In particular, the choice of a $\ring{C}[[\lambda]]$-linear
isomorphism between
$\End_{\boldsymbol{\mathcal{A}}}(\mathcal{E}[[\lambda]],\bullet)$ and
$\End_{\mathcal{A}}(\mathcal{E})[[\lambda]]$ induces a new
\emph{deformed} multiplication $\star'$ for
$\End_{\mathcal{A}}(\mathcal{E})[[\lambda]]$ together with a new
module multiplication $\bullet'$ for $\mathcal{E}[[\lambda]]$
such that $(\mathcal{E}[[\lambda]], \bullet', \bullet)$ becomes a
$(\End_{\mathcal{A}}(\mathcal{E})[[\lambda]],
\star')$-$\boldsymbol{\mathcal{A}}$ bimodule.
\begin{remark}
    \label{remark:CoolStuff}
    \begin{enumerate}
    \item Since $(\mathcal{E}[[\lambda]], \bullet)$ is unique up to
        equivalence the deformation $\star'$ is unique up to
        \emph{isomorphism} since
        $\End_{\boldsymbol{\mathcal{A}}}(\mathcal{E}[[\lambda]],
        \bullet)$ is fixed. One can even obtain a $\star'$ which is
        unique up to equivalence if one imposes $\bullet'$ to be a
        \emph{deformation} of the original left module structure.

        Otherwise, if $\star'$ is such a deformation and $\Phi$ is an
        automorphism of the undeformed algebra
        $\End_{\mathcal{A}}(\mathcal{E})$ then 
        \begin{equation}
            \label{eq:AtwistedB}
            A \star^{\Phi} B := \Phi( \Phi^{-1} A \star' \Phi^{-1}B)
        \end{equation}
        yields another isomorphic but not necessarily equivalent
        deformation of $\End_{\mathcal{A}}(\mathcal{E})$ allowing for
        a bimodule structure as above.
    \item In general, there is an obstruction on $\star'$ to allow
        such a bimodule deformation $\bullet'$ for a given fixed
        $\star$ (and hence $\bullet$) as the algebra structure has to
        be isomorphic to
        $\End_{\boldsymbol{\mathcal{A}}}(\mathcal{E}[[\lambda]],
        \bullet)$.
    \item By analogous arguments as above one can also show the
        existence and uniqueness up to isometries of deformations of
        Hermitian fiber metrics \cite{bursztyn.waldmann:2000b}.
    \item In physical terms: noncommutative field theories on a
        classical vector bundle always exist and are even uniquely
        determined by the underlying deformation of the space-time, at
        least up to equivalence. Morally, this can be seen as the
        deeper reason for the existence of Seiberg-Witten maps.
    \end{enumerate}
\end{remark}

Let us now mention two other constructions leading to deformed vector
bundles. It is clear that the above argument has strong algebraic
power but is of little use when one wants more explicit formulas as
even the classical projections $P_0$ describing a given vector bundle
$E \to M$ are typically rather in-explicit. The following two
constructions provide more explicit formulas:
\begin{enumerate}
\item Jur\v{c}o, Schupp, and Wess \cite{jurco.schupp.wess:2002a}
    considered the case of a \emph{line bundle} $L \to M$ with
    connection $\nabla^L$ and an arbitrary Poisson structure $\theta$
    on $M$. Here one can use first the Kontsevich star product
    quantizing $\theta$ by use of a global formality map.  Second, one
    can use the same formality map together with the connection
    $\nabla^L$ to construct $\bullet$, $\star'$ and $\bullet'$ as
    well. The construction depends on the choice of a global
    formality. One also obtains a Seiberg-Witten map using the
    formality.
\item In \cite{waldmann:2002b} we considered the case of a symplectic
    manifold with arbitrary vector bundle $E \to M$. Given a
    symplectic connection $\nabla$, Fedosov's construction yields a
    star product $\star$ for $M$, see e.g.~\cite{fedosov:1996a}. Using
    a connection $\nabla^E$ for $E$ one obtains $\bullet$, $\star'$
    and $\bullet'$ depending even functorially on the inital data of
    the connections. Hence one obtains a very explicit and geometric
    construction this way.
\end{enumerate}

For both approaches one can show that the resulting deformations
$\bullet$, $\star'$ and $\bullet'$ can be chosen to be \emph{local},
i.e. the deformations are formal power series in bidifferential
operators acting on functions, sections and endomorphisms,
respectively. Thus one can `localize' and restrict to open subsets $U
\subseteq M$. If in particular one has a good open cover
$\{U_\alpha\}$ of $M$ then $E\big|_{U_\alpha}$ becomes a
\emph{trivial} vector bundle. Since the deformation is unique up to
equivalence the restricted deformation $\bullet_\alpha$ has to be
equivalent to the trivial deformation of a trivial bundle. This way
one arrives at a description of $\bullet$, $\bullet'$ and $\star'$ in
terms of \emph{transition matrices} $\boldsymbol{\Phi}_{\alpha\beta}$
satisfying a \emph{deformed cocycle identity}
\begin{equation}
    \label{eq:defcocycle}
    \boldsymbol{\Phi}_{\alpha\beta} \star
    \boldsymbol{\Phi}_{\beta\gamma} \star
    \boldsymbol{\Phi}_{\gamma\alpha} = \Unit
    \quad
    \textrm{and}
    \quad
    \boldsymbol{\Phi}_{\alpha\beta} \star
    \boldsymbol{\Phi}_{\beta\alpha} = \Unit
\end{equation}
on non-trivial overlaps of $U_\alpha$, $U_\beta$, and $U_\gamma$.
Here $\boldsymbol{\Phi}_{\alpha\beta} = \sum_{r=0}^\infty \lambda^r
\boldsymbol{\Phi}^{(r)}_{\alpha\beta} \in
M_k(C^\infty(U_{\alpha\beta}))[[\lambda]]$ and the
$\boldsymbol{\Phi}^{(0)}_{\alpha\beta}$ are the classical transition
matrices.  Conversely, if one finds a deformation
\eqref{eq:defcocycle} of the classical cocycle then one can construct
a deformation of the vector bundle out of it. This can be seen as a
\emph{Quantum Serre-Swan Theorem}, see \cite{waldmann:2002a}.  We
conclude this section with a few further remarks:
\begin{remark}
    \label{remark:KandIndex}
    \begin{enumerate}
    \item Since the finitely generated projective modules
        $\mathcal{E}$ over $\mathcal{A}$ give the $K_0$-theory of the
        algebra $\mathcal{A}$ and since any such $\mathcal{E}$ can be
        deformed in a unique way up to equivalence and since clearly
        any finitely generated projective module over
        $\boldsymbol{\mathcal{A}}$ arises this way up to isomorphism
        one finally obtains that the classical limit map $\cl$ induces
        an \emph{isomorphism}
        \begin{equation}
            \label{eq:KtheoryRigid}
            \cl_*: K_0(\boldsymbol{\mathcal{A}}) 
            \stackrel{\cong}{\longrightarrow} 
            K_0(\mathcal{A}).
        \end{equation}
        Thus $K$-theory is stable under formal deformations
        \cite{rosenberg:1996a:pre}.
    \item If $\int: \boldsymbol{\mathcal{A}} \to \ring{C}[[\lambda]]$
        is a trace functional, i.e.
        \begin{equation}
            \label{eq:TraceFun}
            \int a \star b = \int b \star a,
        \end{equation}
        then $\mathrm{ind}: K_0(\boldsymbol{\mathcal{A}})
        \longrightarrow \ring{C}[[\lambda]]$ defined by
        \begin{equation}
            \label{eq:Index}
            [\boldsymbol{P}] \mapsto \int \tr(\boldsymbol{P})
        \end{equation}
        gives a well-defined group morphism and for a fixed choice of
        $\int$ the \emph{index} $\mathrm{ind}(\boldsymbol{P})$ depends
        only on the classical class $[P_0]$. In case of deformation
        quantization this yields the \emph{index theorems of
          deformation quantization} where one
        has explicit formulas for $\mathrm{ind}(\boldsymbol{P})$ in
        terms of geometric data of $E$, $M$ and the equivalence class
        $[\star]$ of the star product, see Fedosov's book for the
        symplectic case \cite{fedosov:1996a} as well as Nest and
        Tsygan \cite{nest.tsygan:1995a,nest.tsygan:1995b} and the work
        of Tamarkin and Tsygan for the Poisson case
        \cite{tamarkin.tsygan:2001a}.
    \item In the connected symplectic case the trace functional $\int$
        is unique up to normalization \cite{nest.tsygan:1995a} and
        given by a deformation of the integration over $M$ with
        respect to the Liouville measure. In the Poisson case one may
        have many different trace functionals, see e.g.
        \cite{bieliavsky.bordemann.gutt.waldmann:2002a:pre}.
    \item Physically, such trace functionals are needed for the
        formulation of gauge invariant action functionals which are
        used to define dynamics for the noncommutative field
        theories. Recall that the structure of a deformed vector
        bundle is only the kinematical framework.
    \end{enumerate}
\end{remark}

%
%

\section{Morita Equivalence}
\label{sec:morita}

Let us now discuss how Morita theory enters the picture of deformed
vector bundles.  Vector bundles do not only correspond to projective
modules but the projections $P_0 \in M_N(C^\infty(M))$ are always
\emph{full projections} which means that the ideal in $C^\infty(M)$
generated by the components $(P_0)_{ij}$ is the whole algebra
$C^\infty(M)$. We exclude the trivial case $P_0 = 0$ from our
discussion in order to avoid trivialities.  Then the following
statement is implied by general Morita theory, see e.g.
\cite{lam:1999a} as well as \cite{bursztyn.waldmann:2000b}.
\begin{theorem}
    \label{theorem:Morita}
    The bimodule $\mathcal{E} = \Gamma^\infty(E)$ is actually a Morita
    equivalence bimodule for the algebras $C^\infty(M)$ and
    $\Gamma^\infty(\End(E))$. In particular, these algebras are Morita
    equivalent.
\end{theorem}

In the general algebraic case, it is easy to check that the
deformation $\boldsymbol{P}$ of a full projection $P_0$ is again full
whence we conclude that $(\mathcal{E}[[\lambda]], \bullet, \bullet')$
is a Morita equivalence bimodule for the algebras
$(\mathcal{A}[[\lambda]],\star)$ and
$(\End_{\mathcal{A}}(\mathcal{E})[[\lambda]], \star')$ and the later
two algebras are Morita equivalent. Moreover, any Morita equivalence
bimodule between the deformed algebras arises as such a deformation of
a classical Morita equivalence bimodule up to isomorphism, see
e.g.~\cite{bursztyn.waldmann:2002a:pre} for a detailded discussion.
Since the deformation $\star'$ was already fixed up to isomorphism by
the classical right module structure of $\mathcal{E}$, one has to
expect obstructions that an a priori given deformation $\tilde{\star}$
of $\End_{\mathcal{A}}(\mathcal{E})$ is Morita equivalent to the
deformation $\star$ of $\mathcal{A}$. These obstructions make the
classification of the Morita equivalent deformations difficult in the
general framework. We shall come back to this effect when considering
the Picard groupoid.

However, for symplectic star products one has the following explicit
classification of Morita equivalent star products
\cite{bursztyn.waldmann:2002a}, see also
\cite{jurco.schupp.wess:2002a} for a related statement in the Poisson
case. Note that for star products $\star$ and $\star'$ we want the
endomorphisms $\Gamma^\infty(\End(E))$ cassically to be isomorphic to
the functions $C^\infty(M)$ whence the Morita equivalence bimodules
arise as deformations of \emph{line bundles}.
\begin{theorem}
    \label{theorem:ClassSymplME}
    Let $(M, \omega)$ be symplectic. Then two star products $\star$
    and $\star'$ are Morita equivalent if and only if there exists a
    symplectic diffeomorphism $\psi: M \to M$ such that 
    \begin{equation}
        \label{eq:MECond}
        \psi^*c(\star') - c(\star) \in 
        2 \pi\im H^2_{\mathrm{dR}} (M \mathbb{Z}),
    \end{equation}
    where $c(\star) \in \frac{[\omega]}{\im\lambda} +
    H^2_{\mathrm{dR}}(M, \mathbb{C})[[\lambda]]$ is the characteristic
    class of $\star$.  The equivalence bimodule can be obtained by
    deforming a line bundle $L \to M$ whose Chern class $c_1(L)$ is
    given by the above integer class.
\end{theorem}
The most suitable definition of the characteristic class of a
symplectic star product which is used in this theorem is the \v{C}ech
cohomological description as it can be found in
\cite{gutt.rawnsley:1999a}. Then the first proof in
\cite{bursztyn.waldmann:2002a} consists in examining the deformed
transition functions \eqref{eq:defcocycle}. In the approach of
\cite{waldmann:2002b} using Fedosov's construction there is an almost
trivial proof for the above theorem as the Chern class of the line
bundle $L$ can be build into the Fedosov construction as a curvature
term of a connection $\nabla^L$ on $L$ directly.

The additional diffeomorphism $\psi$ is necessary as $\star'$ is only
determined by $L$ up to isomorphism and not up to equivalence as this
is encoded in the characteristic class.

\begin{remark}
    \label{remark:StrongME}
    There is even a stronger result: For $^*$-algebras one has a
    notion of \emph{strong Morita equivalence}
    \cite{bursztyn.waldmann:2001a} which is a generalization of
    Rieffel's notion of strong Morita equivalence for $C^*$-algebras
    \cite{rieffel:1974b}. Applying this for Hermitian star products,
    i.e. those with $\cc{f \star g} = \cc{g} \star \cc{f}$, one has
    the statement that two Hermitian star products are strongly Morita
    equivalent if and only if they are Morita equivalent
    \cite[Thm.~2]{bursztyn.waldmann:2002a}. One uses a deformed
    Hermitian fiber metric in order to get this stronger result.
    Physically, this is the relevant notion of Morita equivalence as
    one also needs to keep track of the $^*$-involutions and
    positivity requirements as we have discussed above. Thanks to
    \cite[Thm.~2]{bursztyn.waldmann:2002a}, we can focus on the purely
    ring-theoretical Morita theory without restriction.
\end{remark}

%
%

\section{The Picard Groupoid}
\label{sec:picard}

In this last section we shall consider the question in `how many ways'
two Morita equivalent algebras can actually be Morita equivalent. In
particular, we want to investigate how Morita equivalence bimodules
behave under formal deformations.

First we note that this is physically an important questions since we
have already seen that the algebra $C^\infty(M)$ and the algebra
$\Gamma^\infty(\End(E))$, which encodes the gauge transformations, are
Morita equivalent via the sections $\Gamma^\infty(E)$ of the vector
bundle $E$. Thus the above question wants to answer how many
`different' vector bundles, i.e. field theories, one can find which
allow for such a bimodule structure for the same algebra of gauge
transformations.

To formulate these questions one uses the following definitions for
unital algebras $\mathcal{A}$, $\mathcal{B}$, \ldots over some ring
$\ring{C}$:
\begin{definition}
    \label{definition:PicC}
    Let $\Picc(\mathcal{A}, \mathcal{B})$ denote the category of
    $\mathcal{A}$-$\mathcal{B}$ Morita equivalence bimodules with
    bimodule homomorphisms as morphisms.  The set of isomorphism
    classes of bimodules in $\Picc(\mathcal{A}, \mathcal{B})$ is
    denoted by $\Pic(\mathcal{A}, \mathcal{B})$.
\end{definition}
From Morita theory we know that $\mathcal{E} \in \Picc(\mathcal{A},
\mathcal{B})$ is a finitely generated and projective module whence the
isomorphism classes are a set indeed.

It is a well-known fact that tensoring equivalence bimodules gives
again an equivalence bimodule. Hence if $\mathcal{E} \in
\Picc(\mathcal{A}, \mathcal{B})$ and $\mathcal{F} \in
\Picc(\mathcal{B}, \mathcal{C})$ then $\mathcal{E}
\otimes_{\mathcal{B}} \mathcal{F} \in \Picc(\mathcal{A},
\mathcal{C})$. Moreover, it is clear that this tensor product is
compatible with the notion of isomorphisms of equivalence
bimodules. Thus this gives a composition law
\begin{equation}
    \label{eq:Tensor}
    \tensor: \Pic(\mathcal{A}, \mathcal{B}) \times \Pic(\mathcal{B},
    \mathcal{C}) \longrightarrow \Pic(\mathcal{A}, \mathcal{C}).
\end{equation}
Then the tensor product is associative on the level of isomorphism
classes, whenever the composition is defined.  We also note that the
trivial self-equivalence bimodule $\mathcal{A}$ behaves like a unit
with respect to $\tensor$, at least on the level of isomorphism
classes. Finally, the dual module to $\mathcal{E}$ gives an inverse
whence we eventually end up with a \emph{groupoid structure}, called
the \emph{Picard groupoid} $\Pic(\cdot,\cdot)$.  The units are just
trivial self-equivalence bimodules and the spaces of arrows are just
the $\Pic(\mathcal{A}, \mathcal{B})$.  The isotropy groups of this
groupoid are the \emph{Picard groups} $\Pic(\mathcal{A}) =
\Pic(\mathcal{A},\mathcal{A})$, see e.g.~\cite{lam:1999a,bass:1968a}.

After this excursion let us now focus again on the deformation
problem. Assume $\boldsymbol{\mathcal{A}} = (\mathcal{A}[[\lambda]],
\star)$ and $\boldsymbol{\mathcal{B}} = (\mathcal{B}[[\lambda]],
\star')$ are associative deformations such that the resulting algebras
are Morita equivalent. Then $\Picc(\boldsymbol{\mathcal{A}},
\boldsymbol{\mathcal{B}})$ is non-empty and any
$\boldsymbol{\mathcal{E}}$ is isomorphic to a deformation
$(\mathcal{E}[[\lambda]], \bullet, \bullet')$ of an equivalence
bimodule $\mathcal{E}$ of the undeformed algebras $\mathcal{A}$,
$\mathcal{B}$. In particular, $\mathcal{A}$ and $\mathcal{B}$ have to
be Morita equivalent, too. Moreover, $\mathcal{E}$ is uniquely
determined up to isomorphism whence one eventually obtains a
well-defined \emph{classical limit map}
\begin{equation}
    \label{eq:clPic}
    \cl_{*}: \Pic(\boldsymbol{\mathcal{A}}, \boldsymbol{\mathcal{B}})
    \longrightarrow \Pic(\mathcal{A}, \mathcal{B}).
\end{equation}
It is easy to see that the classical limit map behaves well with
respect to tensor products of bimodules whence on the level of
isomorphism classes we obtain a groupoid morphism, see
\cite{bursztyn.waldmann:2002a:pre} where the case of the group
morphism is discussed:
\begin{proposition}
    \label{proposition:PicMorph}
    The classical limit map $\cl_{*}$ is a groupoid morphism. In
    particular,
    \begin{equation}
        \label{eq:clPicPic}
        \cl_{*}: \Pic(\boldsymbol{\mathcal{A}}) 
        \longrightarrow \Pic(\mathcal{A})
    \end{equation}
    is a group morphism.
\end{proposition}

Note that this is a very similar situation as for the $K$-theory
\eqref{eq:KtheoryRigid}. However, here $\cl_*$ is far from being an
isomorphism in general. Thus we would like to find a description of
the kernel and the image of the map $\cl_*$, at least for the cases
where $\mathcal{A}$ is commutative.

For the kernel one obtains the following characterization. Let
\begin{equation}
    \label{eq:EquivA}
    \mathrm{Equiv}(\boldsymbol{\mathcal{A}}) = 
    \left\{
        T = \id + \sum_{r=1}^\infty \lambda T_r 
        \; | \; T \in \mathrm{Aut}(\boldsymbol{\mathcal{A}})
    \right\}
\end{equation}
denote the \emph{self-equivalences} of the deformed algebra. Since we
assume that the undeformed algebra $\mathcal{A}$ is commutative, the
inner automorphisms of $\boldsymbol{\mathcal{A}}$ are necessarily
self-equivalences. Thus one can define the group of outer
self-equivalences
\begin{equation}
    \label{eq:OutEquiv}
    \mathrm{OutEquiv}(\boldsymbol{\mathcal{A}}) = 
    \frac{\mathrm{Equiv}(\boldsymbol{\mathcal{A}})}
    {\mathrm{InnAut}(\boldsymbol{\mathcal{A}})}.
\end{equation}
Then one has
\begin{equation}
    \label{eq:kerclOUtEquiv}
    \ker \cl_* \cong \mathrm{OutEquiv}(\boldsymbol{\mathcal{A}})
\end{equation}
as groups \cite[Cor.~3.11]{bursztyn.waldmann:2002a:pre}.

In the case of star products one can describe $\ker\cl_*$ even more
explicitly. Assume that $\star$ is a star product on $(M, \pi)$ with
the property that any $\pi$-central function can be deformed into a
$\star$-central function and any $\pi$-derivation can be deformed into
a $\star$-derivation. There are many star products which actually have
this property, e.g. all symplectic star product, the Kontsevich star
product for a formal Poisson structure which is equal to the classical
one and the star products constructed in
\cite{cattaneo.felder.tomassini:2002a,cattaneo.felder.tomassini:2002b}.
Under these assumptions one has
\cite[Thm.~7.1]{bursztyn.waldmann:2002a:pre}
\begin{equation}
    \label{eq:OutEquivStar}
    \mathrm{OutEquiv}(\star) \cong 
    \frac{H^1_\pi(M, \mathbb{C})}{2 \pi\im H^1_\pi(M, \mathbb{Z})}
    + \lambda H^1_\pi(M, \mathbb{C})[[\lambda]]
\end{equation}
as \emph{sets}, where $H^1_\pi(M, \mathbb{C})$ denotes the first
complex Poisson cohomology of $(M, \pi)$ and $H^1_\pi(M, \mathbb{Z})$
the first integral Poisson cohomology, i.e. the image of the integral
deRham classes under the natural map $H^1_{\mathrm{dR}}(M, \mathbb{Z})
\longrightarrow H^1_\pi(M, \mathbb{C})$.

The identification above is even a \emph{group isomorphism} for
symplectic star products where the right hand side is endowed with its
canonical abelian group structure. However, in the general Poisson
case the group structure on the left hand side is nonabelian.

The situation for the image of the classical limit map $\cl_*$ is more
mysterious \cite{bursztyn.waldmann:2002a:pre}: From the condition
\eqref{eq:MECond} one obtains that the torsion line bundles are always
in the image in the case of symplectic star products. However, there
are examples where the image contains also non-torsion elements and it
seems to depend strongly on the example how big the image actually can
be. In the Poisson case even less is known.

\section*{Acknowledgements}

I would like to thank the organizers of the workshop for their
excellent working conditions and warm hospitality. Moreover, I would
like to thank the participants for their remarks, ideas and
suggestions on the topic as well as Henrique Bursztyn and Stefan
Jansen for comments on the manuscript.

%
%

\begin{footnotesize}

\end{footnotesize}

\end{document}